\documentclass[11pt]{article}
\usepackage{graphicx,amsmath,amssymb,mathrsfs}
\usepackage[a4paper,  margin=1in]{geometry}
\usepackage[utf8]{inputenc}
\usepackage[english]{babel}
\usepackage[utf8]{inputenc}
\usepackage{amsmath}
\usepackage{amsfonts}
\usepackage{amssymb} 
\usepackage{authblk}
\usepackage{natbib}
\usepackage[none]{hyphenat}
\usepackage{setspace}
\setcitestyle{authoryear,open={(},close={)}}

\newtheorem{theorem}{Theorem}[section]
\newtheorem{corollary}{Corollary}[theorem]
\newtheorem{lemma}[theorem]{Lemma}

\begin{document}

\title{\textbf{Bound on FWER for correlated normal distribution  }}

\author[1]{Nabaneet Das}
\author[2]{Subir K. Bhandari}
\affil[1]{Indian Statistical Institute,Kolkata}
\affil[2]{Indian Statistical Institute,Kolkata}

\maketitle

\begin{abstract}
In this paper, our main focus is to obtain an asymptotic bound on the family wise error rate (FWER) for Bonferroni-type procedure in the simultaneous hypotheses testing problem when the observations corresponding to individual hypothesis are correlated. In particular, we have considered the sequence of null hypotheses $H_{0i} : X_i \sim N(0,1) \: (i=1,2,....,n)$ and equicorrelated structure of the sequence $ (X_1,....,X_n)$. Distribution free bound on FWER under equicorrelated setup can be found in \cite{tong2014probability}. But the upper bound provided in \cite{tong2014probability} is not a bounded quantity as the no. of hypotheses ($n$) gets larger and larger and as a result, FWER is highly overestimated for the choice of a particular distribution (e.g.- normal). In the equicorrelated normal setup, we have shown that FWER asymptotically is a convex function (as a function of correlation ($\rho$)) and hence an upper bound on the FWER of Bonferroni-$\alpha$ procedure is $\alpha(1-\rho)$. \\ 
This implies that Bonferroni's method actually controls the FWER at a much smaller level than the desired level of significance under the positively correlated case and necessitates a correlation correction.  
\end{abstract}

\section{ Introduction } 

Multiple hypothesis testing has been one of the most lively area of research in statistics for the past few decades. The biggest challenge in this area comes from the fact that the models involve an extensive collection of unknown parameters and one has to draw simultaneous inference on a large number of hypotheses mainly with the goal of ensuring a good overall performance (rather than focusing too much on the individual problems). Very often data sets from modern scientific investigations in the field of Biology, astronomy, economics etc. require such simultaneous testing on thousands of hypotheses. 
\\
Various measures of error rate have been proposed over the years. One of the hard-line frequentist approach is to control the family wise error rate (FWER) which is defined as the probability of making at least one false rejection in a family of hypothesis-testing problem. \\

Bonferroni's bound provides the classical FWER control method. However, the step-up and step-down algorithms by \cite{holm1979simple}, \cite{simes1986improved}, \cite{,hommel1988stagewise}, \cite{hochberg1988sharper} provides improvement over the Bonferroni's method in terms of power. While Holm's procedure provides control over the FWER in general, the other algorithms depend heavily on the independence of the p-values of the individual hypothesis. \cite{dudoit2003multiple}, \cite{dudoit2007multiple}, \cite{efron2012large} provides excellent review of the whole theory. \\
One of the main limitations of these classical methods that control FWER in the strong sense, is their conservative nature which results in lack of power. A substantial improvement in power has been achieved by considering the False discovery rate (FDR) criterion proposed by  \cite{benjamini1995controlling}. See, for example \cite{benjamini1999step}, \cite{yekutieli1999resampling}, \cite{storey2002direct}, \cite{sarkar2007stepup}, \cite{romano2008control} for further details. \cite{sarkar2008methods}, \cite{efron2012large}   provides an excellent account on the literature on FDR. \\ 
However, most of these works have been done in the context of independent observations. Very little literature can be found that covered correlated variables. \cite{sarkar2008methods} reviews FDR control under dependence. \cite{efron2010correlated} clearly shows the effects of correlation on the summary statistics by pointing out that the correlation penalty depends on the root mean square (rms) of correlations. An excellent review of the whole work can be found in \cite{efron2012large}. All these works gives immense light in the context that FWER or FDR should be treated more carefully where correlation is present. 
\\ 
Some distribution free bound on FWER can be found in \cite{tong2014probability} using Chebyshev-type inequalities. But as these inequalities are distribution free, FWER is highly overestimated for choice of particular distributions (e.g.- normal). Also, these inequalities are not of much use for a large number of hypotheses. \\ 
In our work, we have considered equicorrelated normal distribution and obtained a sharper bound on the FWER for the Bonferroni type FWER control procedure. This distribution also has an important application in modelling the lifetime of coherent systems (see \cite{loperfido2007some}). We have shown that, asymptotically (For large no. of hypotheses)  $ FWER ( \rho) $ is a convex function in $ \rho \in [0,1] $ and hence FWER in Bonferroni-$\alpha $  procedure is bounded by $ \alpha ( 1- \rho) $. This suggests a necessary correlation correction in Bonferroni  procedure. While the bound provided in \cite{tong2014probability} is not a bounded quantity as no. of hypotheses gets larger and larger , the bound provided in this work remains stable even as $ n \to \infty $  and shows a clearer picture of the effect of correlation on FWER. This is probably the first attempt in the context of finding most classically used FWER in terms of $ \rho $ asymptotically. It is also important to mention that the lifetime of parallel systems can be conveniently modelled by the maximum of exchangeable normal random vector. The asymptotic bound on FWER provided in this paper is actually a lower bound on the c.d.f. of the failure time of the parallel systems. \cite{loperfido2007some} has shown that the distribution of the maximum of $n$ observations from exchangeable normal distribution follows $(n-1)$ dimensional skew normal distribution. While the c.d.f. of multivariate skew normal distribution is very difficult to deal with, an asymptotic bound on the c.d.f. can be obtained by a reasoning similar to the one in this paper.

\section{ Description of the problem } 

Let $ X_1,X_2,...... $ be a sequence of observations and the null hypotheses are  
 
$$H_{0i} : X_i \sim N(0,1) \: \: i=1,2,.... $$  

Here we have considered one sided tests (This means, $H_{0i}$ is rejected for large values of $X_i$ (say $X_i > c$)). A classical measure of the type-I error is FWER, which is the probability of falsely rejecting at least one null hypothesis (Which happens if $ X_i > c $ for some $i$ and the probability is computed under the intersection null hypothesis $ H_0= \bigcap\limits_{i=1}^{n} H_{0i} = \{ X_i \sim N(0,1) \: \: \forall \: i=1,2,...,n \}$).
Then,
\begin{center} 
\textbf{FWER} = P(At least one false rejection)= P($ X_i > c $ for some $ i $ $| $  $ H_0 $ )  
\end{center}
Suppose, $ Corr(X_i,X_j) = \rho \: \: \:  \forall \: i \neq j \: \: \: ( \rho \geq 0 ) $. 
\\
Our goal is to provide an asymptotic bound on FWER in terms of $ \rho $. \\

Let, $ H( \rho ) = 1- $ FWER $ ( \rho ) = P( X_i \leq c, \: \: \forall \: i=1,2,....,n \: |\:  H_0) $ 

\section{ Main theorem } 
\begin{theorem}{}
\label{theorem 1}
Suppose each $ H_{0i} $ is being tested at size $ \alpha_n $. If  $ \:  \lim_{ n \to \infty} n \alpha_n = \alpha \in (0,1) $ then,  
as $ n \to \infty , \:   H''( \rho ) \leq   0 $ and hence $ H( \rho ) $ asymptotically is a concave function in [0,1].  
\end{theorem}

\underline{ Note } :- \begin{enumerate}
\item  For $ \rho = 0 $ (Under independence), we must have, FWER = $ 1 - (1- \alpha_n)^n \approx n \alpha_n $. 
\item For $ \rho =1 $ (When $ X_i = X_j $  a.s. $ \forall \:  i \neq j $), we must have FWER = $ \alpha_n $ (Because one rejection would imply rejection of all null hypotheses).  
\end{enumerate}
Suppose $ y = \mathscr{ L } ( \rho) $ denotes the line which joins $ (1, \alpha_n) $ and $ ( 0, 1- (1- \alpha_n)^n) $. The following corollary describes the asymptotic behaviour of FWER as a function of $ \rho $ .  \\

\begin{corollary}
As $ n \to \infty $, FWER $ ( \rho ) $ is bounded above by the line  $ \mathscr{ L } ( \rho) $.  
\end{corollary}

In this section we are going to provide a proof of this theorem. \\

\subsection{ \underline{An alternate form of $ H ( \rho ) $ and it's derivatives} }

Under the framework described above, we can say that under $ H_0 $, the sequence  $\{  X_n \}_{ n \geq 1}  $   is exchangeable. (i.e. $ ( X_{i_1} ,...,X_{ik} ) \sim N_k ( \: \mathbf{0_k} , \: (1 - \rho) I_k + \rho J_k ) $ (Where $ J_k $ is a $ k * k $ matrix of ones)). \\
Then, $   X_k = \theta + Z_k \: , \: \forall \:  k \geq 1 $. \\
Where $ \theta $ is a mean 0, normal random variable, independent of the sequence $ \{ Z_n \}_{ n \geq 1 } $ and $ Z_i $'s are i.i.d. normal random variables. \\
Since $ Cov ( X_i , X_j) = \rho $, this implies that $
Var( \theta )= \rho  $  \\
$ \Rightarrow \theta \sim N(0,\rho) $ and $ Z_n \sim^{ i.i.d.} N(0, 1- \rho) \: \: \forall \:  n \geq 1  $

Thus,
$$ H( \rho ) = P( \theta + Z_i \leq c \: \:  \forall \: i=1,2,...,n )
 = E_{ \theta } \Big[ \Phi^{n} \Big( \frac{ c- \theta}{ \sqrt{ 1 - \rho}} \Big) \Big] = E \Big[ \Phi^{n} \Big( \frac{ c+  \sqrt{ \rho} Z }{ \sqrt{ 1 - \rho}} \Big) \Big]     $$ . \\ 
(Where Z $ \sim N(0,1) $  and $ \Phi $ is the c.d.f. of N(0,1) distribution) 

 \medskip 

If we define, $ d = \frac{ c + \sqrt{ \rho} Z}{ \sqrt{ 1- \rho}} $,
 then  $ H( \rho) = E[ \Phi^{n} (d) ] $.\\

Now, an application of dominated convergence theorem would yield, 
\begin{equation} 
  H ' ( \rho) = E[ n \Phi^{n-1} (d) \phi(d) G( \rho , Z) ] \: \: \:  
 \end{equation}
 (Where $ G( \rho , Z) = \frac{ \partial d}{\partial \rho} = \frac{1}{2 (1- \rho)^{ \frac{3}{2}}} [c+ \frac{Z}{ \sqrt{ \rho}}] $  and $ \phi(.) $ is the N(0,1) p.d.f. ) \\
 
 And again by D.C.T.,
 
 \begin{equation}
   H'' ( \rho) = E \Big[ \frac{n}{2} \Phi^{n-2} (d) \phi(d) \Big( a G^2 ( \rho , Z) + b G( \rho ,Z) + \frac{ c \Phi(d)}{ 4 \rho (1- \rho)^{ \frac{3}{2}}} \Big)   \Big]
  \end{equation} 
   
 Where, $ a= (n-1) \phi(d) - d \Phi(-d) $ and $ b = \frac{ (4 \rho -1) \Phi(d)}{ 2 \rho ( 1 - \rho)} $. \\

Let's define, $ \alpha_1 = \Phi (-d) $. \\
Observe that,
\begin{center} 
 $ H'' ( \rho) = E \Big[ \frac{n}{2} (1- \alpha_1)^{n-2}\phi(d) \Big( a G^2( \rho , Z) + b G( \rho ,Z) + \frac{ c (1- \alpha_1)}{4 \rho (1- \rho)^{ \frac{3}{2}}} \Big) \Big]$ 
\\
$  = \int\limits_{-\infty}^{ \infty}\frac{n}{2} (1- \alpha_1)^{n-2}\phi(d) \Big[ a G^2( \rho , z) + b G( \rho ,z) + \frac{ c (1- \alpha_1)}{4 \rho (1- \rho)^{ \frac{3}{2}}} \Big] \phi(z) dz   $ 
\end{center} 

\section{ Proof of the main theorem } 
The proof of this theorem involves two steps. \\
\begin{itemize}
\item \textbf{ Step 1 :- }  The second and third term in $ H'' ( \rho ) \to 0 $ as $ n \to \infty $. (Proof is given in appendix (lemma \ref{Lemma I }).

\item \textbf{ Step 2 :- } The first term is asymptotically $ \leq 0$. 

\end{itemize}

Proof of step 1 is given in appendix (lemma \ref{Lemma I }).  We shall proceed with the proof of step 2. Before we proceed with the proof, it is important to observe the behaviour of $c$ as  $n \to \infty$. \\ 
\textbf{\underline{Notation}} :-  $ x_n = \Theta ( y_n) $ if $ \exists$ $ c_1 , c_2 > 0 $ and $ M \in \mathbb{N} $ such that, $ c_1 y_n \leq x_n \leq c_2 y_n \: \: \forall \:  n \geq M $. 

Each $ H_{0i} $ is being tested at size $ \alpha_n $ and we reject $ H_{0i} $ if $ X_i > c $. So, $ \alpha_n = \Phi(-c) $. \\

As $ n \to \infty $ we must have $ c \to \infty $ by the condition $ \lim_{n \to \infty} n \alpha_n = \alpha \in (0,1) $. \\

For large $ c $, we have $ \Phi(-c) \sim \frac{ \phi (c) }{c} = \frac{1}{ \sqrt{2 \pi} c e^{ \frac{c^2}{2}}}. $ \\

So, $ n \sim \alpha \sqrt{2 \pi} c e^{ \frac{c^2}{2}} $  and from this, we can conclude that, $c^2 = \Theta ( \log n ) $. 
\subsection{ \underline{Proof of step 2} }

Partition the range of $\alpha_1$ according as $\{  \alpha_1 \geq \frac{1 }{n} \}  $ and  $\{ \alpha_1 < \frac{1}{n}  \} $. We shall show that the integrals in the first region  $ \to 0 $ as $ n \to \infty$ and the integral in the second region is asymptotically non-positive. 

\begin{itemize}

\item \textbf{ Case 1 } : -  $ \{  \alpha_1 \geq  \frac{1}{n }  \} $ \\
\vspace{0.07 in}

Suppose at $z=z_0, \: \: \alpha_1 (z_0) = \Phi( - d( z_0)) = \frac{1}{n} $.\\
First we'll show that $z_0$ takes a very large negative value. Observe that, $$ \Phi ( - d(z_0) ) = \frac{1}{n} \text{ and } \Phi( - c) = \alpha_n \sim \frac{ \alpha}{n} $$ 
So, $ \: \frac{ \Phi( - d (z_0) ) }{ \Phi ( -c) } \sim \frac{1}{ \alpha} $.  
We know that, $ x \Phi(-x) \sim \phi(x)  $  for large enough $x$.\\ 

So, for large enough $n$,  $\: d(z_0) \Phi( - d(z_0) ) \sim \phi( - d (z_0)  $ and $c \Phi( -c ) \sim \phi(c) $. \\ 
Thus, 
$$ \frac{ \Phi( - d (z_0) ) }{ \Phi ( -c) } \sim \frac{ \frac{ \phi( d(z_0) )}{ d(z_0) } }{ \frac{ \phi (c) }{ c} } = \frac{ c e^{ \frac{c^2}{2} }}{ d(z_0) e^{ \frac{ d^2 (z_0 ) }{ 2} } } \sim \frac{1}{ \alpha} $$  
And this implies, 
$$ \log \Big( \frac{c}{ d (z_0)} \Big)  + \Big( \frac{ c^2 - d^2 (z_0 )}{2} \Big) \sim \log ( \frac{1}{ \alpha} )  $$ 
Since $0 < \alpha < 1 $, we must have $ c > d(z_0) $. Thus,  
$ \frac{c^2 - d^2 (z_0) }{ 2}  \leq \log ( \frac{1}{ \alpha } ) $ \\
\vspace{0.05 in} 
We can now conclude that, $  c - d(z_0 )  \leq \frac{ 2 \log ( \frac{1}{ \alpha} ) }{ c + d(z_0) } $  for large enough $n$. \\ 
Since both $c \to \infty $ and  $d(z_0) \to \infty$ as $n \to \infty$, this implies that $c - d(z_0) \to 0 $ as $n \to \infty$. \\ 

From this, it is easy to see that, $ z_0 + c T( \rho) \to 0 $ as $n \to \infty$. \\ 
(Where $ T( \rho) = \frac{1}{ 1 + \sqrt{ 1 - \rho } }  $ is entirely dependent on $\rho$ and is a positive function for all $0 \leq \rho \leq 1 $.)

This means, $ z_0 $ takes a very large negative value. \\ 

Consider the region $ \{ \alpha_1 \geq \frac{1}{n}  \} $. Since $\alpha_1$ is a decreasing function of $z$, this region can also be written as $\{ z \leq z_0 \} $ \\
 
Note that, $ d G^2 ( \rho , Z) = p (Z,c) $ (A polynomial in $Z$ and $c$). Using this, it is easy to see that, $ \sup\limits_{ z \leq z_0 } d G^2 ( \rho , z) \phi (z) \to 0 $ as $n \to \infty$. 

\vspace{0.03 in}

$d=\frac{c+\sqrt{\rho} Z}{\sqrt{ 1- \rho}} \overset{a.s.}{\to} \infty$ as $ n \to \infty $ (Because $c \to \infty$ as $ n \to \infty$).\\
\vspace{0.03 in} 
So,  $ a = (n-1) \phi(d) - d \Phi(d) \sim d ( n \alpha_1 -1)  $. \\
By lemma \ref{Lemma II} in appendix, we can say that,  
$$ E \Big[\frac{n}{2} (1- \alpha_1)^{n-2} \phi(d)a G^2 ( \rho , Z)\Big] \: \sim \:  E \Big[ \frac{n}{2} (1- \alpha_1)^{n-2} \phi(d) d ( n \alpha_1 -1) G^2 ( \rho , Z) \Big]  $$ 
Now observe that, \\ 

$E \Big[ \frac{n}{2} (1- \alpha_1)^{n-2} \phi(d) d ( n \alpha_1 -1) G^2 ( \rho , Z) \Big] $ \\ 
$ = \int\limits_{ - \infty}^{\infty} n ( n \alpha_1 -1) (1- \alpha_1)^{n-2} \phi(d) d G^2 ( \rho, z) \phi(z) dz $  \\

 $  \propto \int\limits_{0}^{1} f_{n} ( \alpha_1)  K( \alpha_1) d \alpha_1 $  \\ 
 (Where $ f_n ( \alpha_1) = n ( n \alpha_1 -1) (1- \alpha_1)^{n-2} $ and  $ K( \alpha_1) = d G^2 ( \rho , z) \phi(z) $) \\

 Observe that, $ f_n ( \alpha_1) >(<) \: 0  \: \Leftrightarrow \:  \alpha_1 >(<)  \: \frac{1}{n} $ and $ \int\limits_{0}^{1} f_n ( \alpha_1) d \alpha_1 = 0 $. \\ 
 
And, $ \int\limits_{ \frac{1}{n}}^{1} f_n ( \alpha_1)  d \alpha_1 = \int\limits_{1}^{n} (k-1) (1- \frac{k}{n})^{n-2} dk \leq \int\limits_{1}^{\infty} (x-1) e^{ - \frac{x}{2} } dx < \infty $ \\

$ \Rightarrow  \int\limits_{ 0}^{1} |f_n( \alpha_1) | d \alpha_1 $ is bounded.  \\

Since $\sup\limits_{ \alpha_1 \geq \frac{1}{n} } K( \alpha_1) = \sup\limits_{ z \leq z_0 } d G^2 ( \rho ,z) \phi (z) \to  0 $, it is easy to deduce that,
$$ E \Big[ n (1- \alpha_1)^{n-2} \phi(d) d (n \alpha_1 -1) G^2 ( \rho ,Z) I_{ ( \alpha_1 \geq \frac{1}{n} ) } \Big] = \int\limits_{ \frac{1}{n}}^{1} f_n ( \alpha_1) K( \alpha_1) d \alpha_1 \to 0 $$ 
And this implies, 
$$ E \Big[\frac{n}{2} (1- \alpha_1)^{n-2} \phi(d)a G^2 ( \rho , Z) I_{( \alpha_1 \geq \frac{1}{n}) } \Big] \to 0 \text{ as } n \to \infty $$

\item \textbf{ Case 2 } : -  $ \{  \alpha_1 <   \frac{1}{n }  \} $ \\ 

\vspace{ 0.05 in} 

If $\alpha_1 < \frac{1}{n} $ then $ f_n ( \alpha_1) < 0  $. This means $f_{n } ( \alpha_1) K( \alpha_1) < 0 $ and hence 
$$ E \Big[ n (1- \alpha_1)^{n-2} \phi(d) d (n \alpha_1 -1) G^2 ( \rho ,Z) I_{ ( \alpha_1 <  \frac{1}{n} ) } \Big] = \int\limits_{0}^{ \frac{1}{n}} f_n ( \alpha_1) K( \alpha_1) d \alpha_1 \leq 0 $$ 
By lemma \ref{Lemma II}, we can say that, 
$$ E \Big[\frac{n}{2} (1- \alpha_1)^{n-2} \phi(d)a G^2 ( \rho , Z) I_{( \alpha_1 < \frac{1}{n}) } \Big] \leq 0 \text{ for large } n $$

This completes the proof of step 2.

\end{itemize}
This completes the proof of main theorem. 

\section{ Conclusion } 

From the theorem in the previous section, we have $ H''( \rho) \leq 0 $ asymptotically. \\
$ \Rightarrow $  $ FWER''( \rho ) \geq 0 $ as $ n \to \infty. $ Thus, $ FWER ( \rho ) $ asymptotically is a convex function and hence $ FWER  ( \rho) $ is bounded by $ \mathscr{L} ( \rho) $ in [0,1]. From this, we can conclude the following theorem. 
\begin{theorem}{}
For large $n$,  $FWER  ( \rho)  \leq \alpha_n + ( 1- \rho) [ 1- \alpha_n - (1-\alpha_n)^n] $. 
\end{theorem}

For large $n$, $ 1- ( 1- \alpha_n)^n \approx n \alpha_n $ and this implies that FWER $ ( \rho ) \leq \alpha_n [ n - (n-1) \rho] $. \\

Bonferroni's method suggests us to take $ \alpha_n = \frac{ \alpha}{n} $  if we want to maintain $ \alpha $ FWER level. This satisfies the criterion of the main theorem of section 2.  \\

When $ \alpha_n = \frac{ \alpha}{n} $, then $ \alpha_n ( n - (n-1) \rho) \sim \alpha ( 1- \rho) $. \\

Thus, the FWER of Bonferroni's procedure is asymptotically bounded by $ \alpha ( 1- \rho ) $.  

\section{ Appendix }

\begin{lemma}
\label{Lemma I }
The second and third term in $ H''( \rho) \to 0 $ as $ n \to \infty. $
\end{lemma} 
\textbf{ Proof} :-  We shall do it only for the third term. The other one follows similarly.

Third term in $ H''( \rho) $ is proportional to $  E [ nc (1- \alpha_1)^{n-1} \phi(d) ]  = \int\limits_{ - \infty}^{ \infty} nc (1- \alpha_1)^{n-1} \phi(d) \phi(z) dz $.\\

Let, $ \phi (z) = M( \alpha_1) $. \\
\vspace{0.05 in} 

Then, $ \int\limits_{ - \infty}^{ \infty} nc (1- \alpha_1)^{n-1} \phi(d) \phi(z) dz  = \int\limits_{0}^{1} nc (1- \alpha_1)^{n-1} M( \alpha_1) d \alpha_1 $\\

\vspace{0.05 in} 
Since $ \phi(d) \leq d^3 \alpha_1 $ for large enough $d$ and 
$c^2  = \Theta ( \log n) $, we can conclude that, 
$$ E \Big[ n c ( 1- \alpha_1)^{n-1}  \phi(d) I_{ ( \alpha_1 < \frac{1}{ n ( \log n)^3} )} \Big] \to 0 $$ 

Also, if  $ \alpha_1 > \frac{ 6 \log n}{n} $  then $ (1- \alpha_1)^{n-1} \leq \frac{1}{n^3} $. Since $\frac{c}{n^2} \to 0$, we can say that, 
$$  E \Big[ n c ( 1- \alpha_1)^{n-1}  \phi(d) I_{ ( \alpha_1 > \frac{ 6 \log n}{ n } )} \Big] \to 0 $$

Let, $ \Phi ( - d( z_0 - \delta_{n}^{'} ) ) = \frac{1}{ n ( \log n)^3} $ and $ \Phi ( - d( z_0 + \delta_{n}^{''} ) ) = \frac{6 \log n }{ n } $.  \\ 

\vspace{0.05 in} 

The region $\: \{ \frac{1}{ n ( \log n)^3} \leq \alpha_1 \leq \frac{ 6 \log n}{n} \} $ can be written as $ \: z \in [ z_0 - \delta_{n}^{'} , z_0 + \delta_{n}^{''} ]$. \\ 
\vspace{0.05 in} 
An argument similar to the one given in the proof of $c - d(z_0) \to 0 $ in the step 2 of theorem \ref{theorem 1}, it can be deduced that, $ \max \{ \delta_{n}^{'} , \delta_{n}^{''} \} \to 0 $ as $ n \to \infty$. \\

Since $ z_0  + T( \rho ) c \: \to 0 $, for some $T( \rho) > 0$, we can conclude that, 
$$\max\limits_{ z \in [ z_0 - \delta_{n}^{'} , z_0 + \delta_{n}^{''} ]} c  \phi(z) = c \phi ( z_0 + \delta_{n}^{''}) \leq c e^{ - z_{0}^2 }  $$  

Thus, if  $ \: \frac{1}{ n ( \log n)^3} \leq \alpha_1 \leq \frac{ 6 \log n}{n}  $,   then  $ c M ( \alpha_1 )  \leq \max\limits_{ z \in [ z_0 - \delta_{n}^{'} , z_0 + \delta_{n}^{''} ]} c  \phi(z) 
\leq c e^{ - z_{0}^2 }  $\\ 

 Now observe that, \\ 
 
 \vspace{0.05 in} 
 
 $ E \Big[ nc (1- \alpha_1)^{n-1} \phi(d) I_{ \big(\frac{1}{ n (\log n)^3 } \leq  \alpha_1 \leq  \frac{ 6 \log n}{n} \big)} \Big] $ \\ 
 
 $= \int\limits_{\frac{1}{ n  ( \log n)^3 } }^{ \frac{ 6 \log n }{n} } nc  (1- \alpha_1)^{n-1} M( \alpha_1) d \alpha_1  $ \\
 
 $ \leq c e^{ - z_{0}^2} \int\limits_{\frac{1}{ n  ( \log n)^3 } }^{ \frac{ 6 \log n }{n} } n  (1- \alpha_1)^{n-1} d \alpha_1 $ \\ 
 
 $ \leq 6 c\log ( n ) e^{ - z_{0}^2 }  $  \\ 
 
 \vspace{0.05 in} 
 
 We have already shown that, $c^2 = \Theta ( \log n) $ and $ |z_0 | \sim T( \rho) c $. Hence, $6 c \log (n ) e^{ - z_{0}^2 } \to 0 $ and as a direct implication, we have that,  
 $$ E \Big[ nc (1- \alpha_1)^{n-1} \phi(d) I_{ \big( \frac{1}{ n (\log n)^3 } \leq \alpha_1 \leq \frac{ 6 \log n}{n} \big)} \Big]  \to 0  \text{  as  } n \to \infty $$
 
 Hence, third term in $ H'' ( \rho) \to 0 $ as $ n \to \infty$.

\begin{lemma}
\label{Lemma II}
$E \Big[ \frac{n}{2} (1- \alpha_1)^{n-2} \phi(d) G^2 ( \rho , Z) | a - d ( n \alpha_1 -1) |\Big] \to 0 $ as $ n \to \infty$.
\end{lemma}  

\textbf{ Proof} :-  $ | a - d (n \alpha_1 -1) | = (n-1) | \phi(d) - d \Phi(-d) | $. \\
If $ d \leq 2 $, then $ \alpha_1 \geq \Phi (-2) \geq 0.02 $. Hence, $ ( 1- \alpha_1)^{n-2} \leq (0.98 )^{n-2} $. \\
Since $ n = \Theta ( c e^{ \frac{c^2}{2} } ) $, this immediately implies that, 
$$ E \Big[ \frac{n}{2} (1- \alpha_1)^{n-2} \phi(d) G^2 ( \rho , Z) | a - d ( n \alpha_1 -1) | I_{ ( d \leq 2 ) } \Big] \to 0 $$ 

When $ d>2 , \: \:  | \phi(d) - d \Phi(-d) | \leq d \alpha_1  $ and hence $| a - d ( n \alpha_1 -1) | \leq (n-1) d \alpha_1 $. \\

This implies that,
$$ E \Big[ \frac{n}{2} (1- \alpha_1)^{n-2} \phi(d) G^2 ( \rho , Z) | a - d ( n \alpha_1 -1) | I_{ ( d > 2 ) } \Big]  
\leq  E \Big[ n^2 \alpha_1 ( 1- \alpha_1)^{n-2} d \phi(d) G^2 ( \rho , Z) I_{ ( d >2 ) } \Big] $$
An idea similar to the proof of lemma \ref{Lemma I } will tell us that, we need to consider the region $ \{ \frac{1}{n (\log n)^3} \leq \alpha_1 \leq \frac{ 8 \log n}{n} \}$ only. It can be shown in the similar way that, the integral corresponding to this region also $ \to 0 $. \\ 

\section{Simulation results} 
Theorem 5.1 tells us for large $n$, FWER for Bonferroni's method (with level of significance $\alpha$) is asymptotically bounded above by $\alpha (1- \rho)$. In order to verify this result empirically, some simulation results have been provided in table 1. In our simulation experiments, we have considered $ \rho = 0,0.1,0.3,0.5,0.7,0.9$ and $ \alpha=0.01,0.05,0.1,0.4,0.6,0.7$. In each combination of $(\rho,\alpha)$, $10000$ replications have been made to estimate the FWER (the estimate obtained is denoted by $\hat{FWER}$). In each replication, we have generated $10000$ equicorrelated normal random variables each with mean 0 and variance 1. Bonferroni's method suggests us to reject $H_{0i}$ at level $\alpha$ if $Z_i > (1- \frac{ \alpha}{10000})$-th quantile of N(0,1) distribution. In each replication we have to note whether or not any of the $10000$ $Z_i$'s exceeds that cut-off and then  $\hat{FWER}$ is obtained accordingly from the $10000$ replications.\\ 
Each $\hat{ FWER} $ obtained at the combination $(\rho,\alpha)$ is compared with $\alpha(1-\rho)$(the upper bound mentioned in section 5). It is impressive that in all the cases $\hat{FWER}$ is substantially smaller than $\alpha(1-\rho)$ (except at $(\rho,\alpha)=(0.1,0.01)$, although the difference is not noteworthy). All these observations suggest that in positively correlated setup, Bonferroni's method actually controls the FWER at a much smaller level than the desired level of significance which makes this method more conservative.

\begin{table}[h!]
\centering
\begin{tabular}{||c |  c c c c c c c||} 
 \hline
 $\rho $ & $ \alpha $  & 0.01 & 0.05 & 0.1 & 0.4 & 0.6 & 0.7 \\
  \hline
 0.9 & $\hat{FWER} $ & 9.00E-05 & 0.00046  & 0.00053 & 0.00221 & 0.00324 & 0.0031
\\ 

  & $ \alpha (1- \rho) $ & 1.00E-03 & 0.005 
& 0.01
& 0.04
& 0.06
& 0.07

 \\
 \hline
  0.7 & $\hat{FWER}$ & 0.00101 & 0.00363
& 0.00588
& 0.01617
& 0.02149
& 0.023
\\ 

   & $ \alpha (1- \rho) $ & 0.003 & 0.015
& 0.03
& 0.12
& 0.18
& 0.21
  \\
   
\hline
 0.5 & $\hat{FWER}$ & 0.00347 & 0.01156
& 0.01918
& 0.04909
& 0.06414
& 0.07042
\\ 

   & $ \alpha (1- \rho) $  & 0.005 & 0.025
& 0.05
& 0.2
& 0.3
& 0.35
 \\ 
\hline    
 0.3 & $\hat{FWER}$ & 0.00683 & 0.02523
& 0.04363
& 0.11495
& 0.15013
& 0.16494
 \\ 

   & $ \alpha (1- \rho) $ & 0.007 & 0.035
& 0.07
& 0.28
& 0.42
& 0.49
 \\ 
\hline
 0.1 & $\hat{FWER}$ & 0.00996 & 0.04367
& 0.07978
& 0.23801
& 0.31105
& 0.34295
\\ 

   & $ \alpha (1- \rho) $ & 0.009 & 0.045
& 0.09
& 0.36
& 0.54
& 0.63
\\ 
\hline

 0 & $\hat{FWER}$ & 0.01018 & 0.0486 & 0.09424 & 0.32914 & 0.45065 & 0.50499 \\ 

   & $ \alpha (1- \rho) $  & 0.01 & 0.05 & 0.1 & 0.4 & 0.6 & 0.7 \\ 
   
\hline 

\end{tabular}
\caption{Simulation results}
\label{table:1}
\end{table}

\newpage

\nocite{romano2008control,sarkar2002some,efron2007correlation,efron2010correlated,efron2012large,efron2001empirical,efron2007size}

\bibstyle{plain}
\bibliography{references.bib}

\end{document}